\documentclass[11pt, reqno]{amsart}
\usepackage{amsmath}
\usepackage{amssymb}
\usepackage{amsfonts}
\usepackage[usenames]{color}
\usepackage{version}
 \newtheorem{theorem}{Theorem}[section]
 \newtheorem{corollary}[theorem]{Corollary}
 \newtheorem{lemma}[theorem]{Lemma}
 \newtheorem{proposition}[theorem]{Proposition}

\newtheorem{definition}[theorem]{Definition}
\newtheorem{remark}[theorem]{Remark}
\newtheorem{example}[theorem]{Example}
\newtheorem{fact*}{Fact}

\DeclareMathOperator{\diag}{diag}

\DeclareMathOperator\spa{span}

\DeclareMathOperator\Red{Red}
\newcommand\half{\tfrac 12}

\newcommand{\bb}{B^2}

\newcommand{\F}{\mathcal F}

\newcommand{\mm}{\chi} 
\newcommand{\M}{\mathcal{M}}

\newcommand{\vv}{\omega} 

\newcommand{\D}{\mathbb{D}}
\newcommand{\C}{\mathbb{C}}

\newcommand{\ran}[1]{\operatorname{ran}#1}

\newcommand{\ip}[2]{\left\langle #1, #2 \right\rangle}

\newcommand{\inv}{^{-1}}

\newcommand{\p}{\mathcal{P}_2}

\newcommand{\ph}{\varphi}

\newcommand\ga{\gamma}

\newcommand\la{\lambda}

\newcommand{\1}{\mathbf{1}}
\newcommand\beq{\begin{equation}}
\newcommand\ds{\displaystyle}
\newcommand\eeq{\end{equation}}
\newcommand\df{\stackrel{\rm def}{=}}

\newcommand\black{\color{black}}

\newcommand\bbm{\begin{bmatrix}}
\newcommand\ebm{\end{bmatrix}}
\numberwithin{equation}{section}

\begin{document}

\title{Symmetric functions of two noncommuting variables}

\author{J. Agler and N. J. Young}

\address{Department of Mathematics, University of California at San Diego, San Diego, CA \textup{92103}, USA \\
 {\em and} School of Mathematics, Leeds University, Leeds LS2 9JT, U.K.~~
{\em and} School of Mathematics and Statistics, Newcastle University, Newcastle upon Tyne NE1 7RU, U.K.}


\date{5th July 2013}

\begin{abstract}
We prove a noncommutative analogue of the fact that every symmetric analytic function of $(z,w)$ in the bidisc $\D^2$ can be expressed as an analytic function of the variables $z+w$ and $zw$.  We construct an analytic nc-map $S$ from the biball to an infinite-dimensional nc-domain $\Omega$ with the property that, for every bounded symmetric function $\ph$ of two noncommuting variables that is analytic on the biball, there exists a bounded analytic nc-function $\Phi$ on $\Omega$ such that $\ph=\Phi\circ S$.  We also establish a realization formula for $\Phi$, and hence for $\ph$, in terms of operators on Hilbert space.
\end{abstract}
\maketitle

\section{Introduction}\label{intro}
Every symmetric polynomial in two commuting variables $z$ and $w$ can be written as a polynomial in the variables $z+w$ and $zw$; conversely every polynomial in $z+w$ and $zw$ determines a symmetric polynomial in $z$ and $w$.  A similar assertion holds for symmetric analytic functions on symmetric domains in $\C^2$.  For noncommuting variables, on the other hand, no such simple characterizations are valid.  For example, the polynomial
\[
zwz+wzw
\]
in noncommuting variables $z,w$ { \em cannot } be written as $p(z+w,zw+wz)$ for any polynomial $p$;
 M. Wolf showed in 1936 \cite{wolf} that there is no finite basis for the ring of symmetric noncommuting polynomials over $\C$.  She gave noncommutative analogues of the elementary symmetric functions, but they are infinite in number.

In this paper we extend Wolf's results from polynomials to symmetric analytic functions in noncommuting variables within the framework of noncommutative analysis, as developed by J. L. Taylor \cite{taylor} and many other authors, for example \cite{AlpKal,BGM,HelKlepMcC,KalVin,popescu2006,voic}.  We prove noncommutative analogues of the following simple classical result.

Let $\pi:\C^2\to\C^2$ be given by
\[
\pi (z,w)=(z+w,zw).
\]
If $\ph: \D^2\to \C$ is analytic and symmetric in $z$ and $w$ then there exists a unique analytic function $\Phi:\pi(\D^2)\to \C$ such that the following diagram commutes:
\begin{equation*}
\begin{array}{ccccccccc}
  \D^2 \!\!\! & ~ &\stackrel{\pi}{\longrightarrow}& 
~ & \!\!\! \pi(\D^2)\!\!\!  \\ 
~ &  \vcenter{\llap{$\scriptstyle{\ph}$}}\searrow & ~ 
& \swarrow\vcenter{\rlap{$\scriptstyle{\Phi}$}} & ~ \\ 
   ~ & ~ & \!\!\! \C\!\!\!  & ~ & ~ \\  
\end{array}
\end{equation*}
In this diagram the domain $\pi(\D^2)$ is two-dimensional, in consequence of the fact that there is a basis of the ring of symmetric polynomials consisting of two elements, $z+w$ and $zw$.   Wolf's result implies that in any analogous statement for symmetric polynomials in two {\em noncommuting} variables, $\pi(\D^2)$ will have to be replaced by an infinite-dimensional domain.  The same will necessarily be true for the larger class of symmetric holomorphic functions of two noncommuting variables.

We use the notions of {\em nc-functions} and {\em nc-maps} on {\em nc-domains}, briefly explained in Section \ref{ncfunctions}.  An example of an nc-domain is the {\em biball} 
\[
\bb\stackrel{\rm def}{=} \bigcup_{n=1}^\infty {B}_n \times {B}_n,
\]
where ${B}_n $ denotes the open unit ball of the space $\mathcal{M}_n$  of $n\times n$ complex matrices.   
$\bb$ is the noncommutative analogue of the bidisc.
Another example of an nc-domain is the space
\[
\M^\infty \stackrel{\rm def}{=} \bigcup_{n=1}^\infty \M_n^\infty
\]
of infinite sequences of $n\times n$ matrices, for any $n\geq 1$.

The following result is contained in Theorem \ref{realizn} below.
\begin{theorem}\label{realizn1}
  There exists an nc-domain $\Omega$ in $\M^\infty$ such that the map $S:\bb\to \M^\infty$ defined by
\beq\label{defS1}
S( x)= (u,v^2,vuv,vu^2v,\dots),
\eeq
where
\beq\label{defuv1}
u=\frac{x^1+x^2}{2}, \qquad v=\frac {x^1-x^2}{2},
\eeq
has the following two properties.
\begin{enumerate}
\item $S$  is an analytic nc-map from $\bb$ to $\Omega$;
 \item for every bounded symmetric nc-function $\ph$ on the biball there exists a bounded analytic nc-function $\Phi$ on $\Omega$ such that the following diagram commutes:
\begin{equation*}
\begin{array}{ccccccccc}
  B^2 \!\!\! & ~ &\stackrel{S}{\longrightarrow}& 
~ & \!\!\! \Omega\!\!\!  \\ 
~ &  \vcenter{\llap{$\scriptstyle{\ph}$}}\searrow & ~ 
& \swarrow\vcenter{\rlap{$\scriptstyle{\Phi}$}} & ~ \\ 
   ~ & ~ & \!\!\! \bigcup_n \mathcal{M}_n\!\!\!  & ~ & ~ \\  
\end{array}
\end{equation*}
\end{enumerate}

Moreover $\Phi$ can be expressed by the formula
\[
\Phi = \mathcal{F} \circ \Theta_U
\]
for some graded linear fractional transformation $\mathcal{F}$ and some unitary operator $U$ on $\ell^2$, where $\Theta_U$ denotes the functional calculus corresponding to $U$.
\end{theorem}
The sense in which the maps $\ph, \ S$ and $\Phi$ are analytic is explained in Definitions \ref{2.1}, \ref{ncdomaininfty} and  \ref{defncfunc} in the next section; graded linear fractional transformations are explained in Section \ref{linfrac}.

  The domain $\Omega$ of Theorem \ref{realizn1} is not the analogue of the symmetrized bidisc in all respects:  $S$ is far from surjective onto $\Omega$, and we make no uniqueness statement for $\Phi$ in the theorem. 

A more algebraic approach to symmetric functions in noncommuting variables has been adopted by many authors, for example, I. M. Gelfand {\em et al.} \cite{gelfand}.  In the latter paper the action of the symmetric group of order two on polynomials differs from the action studied in the present paper (see \cite[Example 7.16]{gelfand}).

\section{nc-Functions}\label{ncfunctions}
The settings for nc-functions are the ``universal spaces" $\M^d$ comprising $d$-tuples of matrices of all orders, where $d$ is a positive integer or $\infty$.  For $n$ in the set $\mathbb{N}$ of natural numbers
we denote by $\M_n$ the space of $n\times n$ complex matrices with the usual operator norm.  For $1\leq d<\infty$ the space $\M_n^d$ of $d$-tuples of $n\times n$ matrices is a Banach space with norm
\[
\| (M^1,\dots,M^d)\| = \max_{j=1,\dots,d} \|M_j\|.
\]
For $d=\infty$ it is more convenient to index sequences by the non-negative integers, so that a typical element of $\M_n^\infty$ will be written $g=(g^0,g^1,g^2,\dots)$ with $g^j\in \M_n$.  Of course $\M_n^\infty$ is not naturally a normed space, but it is a Fr\'echet space with respect to the product topology.

For $d\leq \infty$ define 
\[
\M^d \df \bigcup_{n=1}^\infty \M_n^d.
\]
A set $U\subset \M^d$ is said to be 
 {\em nc-open} if $U\cap \M^d_n$ is open in $\M^d_n$ for every $n\geq 1$.  When $d<\infty$ the space $\M^d$ is a disjoint union of Banach spaces.

\begin{definition}\label{2.1}
Let $d\in\mathbb{N}$.  An {\em nc-domain} in $\M^d$ is a subset $D$ of $ \M^d$ that is nc-open and satisfies
\begin{enumerate}\label{ncdomain}
\item if $M, N \in D$ then $M\oplus N\in D$, and
\item if $M\in D\cap \M_m^d$ and $U\in \M_m$ is unitary then $U^*MU \in D$.
\end{enumerate}
\end{definition}
Here if $M=(M^1,\dots,M^d)\in \M_m^d$ and $N=(N^1,\dots,N^d) \in \M_n^d$
then $M\oplus N$ denotes $(M^1\oplus N^1,\dots,M^d\oplus N^d)\in \M_{m+n}^d$, where $M^j \oplus N^j$ is the $(m+n)$-square block diagonal matrix $\diag(M^j,N^j)$.  In (2) $U^*MU$ denotes $(U^*M^1U, \dots,U^*M^dU)$.

The nc-domains are the natural domains on which to define nc-functions -- see below.  For $d=\infty$ it is too restrictive to require that nc-domains be nc-open: there are too few nc-open sets.  The following refinement is a more fruitful notion.
\begin{definition}\label{ncdomaininfty}
An {\em nc-domain} in $\M^\infty$ is a subset $D$ of $ \M^\infty$ that is open in some 
 union of Banach spaces contained in $\M^\infty$ and satisfies conditions {\rm (1)} and {\rm(2)} of Definition {\rm \ref{ncdomain}}.

Here a {\em union of Banach spaces contained in $\M^\infty$}  is a subset $E$ of $\M^\infty$ such that, for every $n\in\mathbb{N}$,  $E\cap \M_n^\infty$ is a Banach space with respect to some unitarily invariant norm $\|\cdot\|_n$ that induces a finer topology than the product topology on $E\cap \M_n^\infty$, and $E$ carries the topology of the disjoint union of the spaces $(E\cap \M_n^\infty, \|\cdot\|_n)_{n\geq 1}$.
\end{definition}
Here of course to say that $\|\cdot\|_n$ is unitarily invariant on $E\cap \M_n^\infty$ means that, for every $x\in E\cap \M_n^\infty$ and every invertible matrix $u\in\M_n$,
\[
\|u^* xu\|_n = \|x\|_n.
\]

\begin{example}\label{exAofD} \rm The {\em nc-disc algebra} $\mathbf{A}(\D)$ is the space of analytic square-matrix-valued functions on $\D$ that extend continuously to the closure of $\D$, with the supremum norm.    The space $ \mathbf{A}(\D)$ is a union of Banach spaces  contained in $\M^\infty$ (see Proposition \ref{boldAofD} below).
\end{example}
\begin{definition}\label{defncfunc}
An {\em nc-function} is a function $\ph: D\to \M^1$ for some nc-domain $D$ in $\M^d$ (for some $d\leq\infty$) that satisfies the conditions
\begin{enumerate}
\item $\ph$ maps $D\cap\M_n^d$ to $\M_n$ for every $n\in\mathbb{N}$;
\item for all $M,N \in D $, 
\beq\label{oadditive}
\ph(M\oplus N) = \ph(M)\oplus \ph(N),
\eeq
 and
\item for all $n\in\mathbb{N}$, all $M\in D\cap \M_n^d$ and all invertible matrices $s\in\M_n$ such that $s\inv M s\in D$,
\beq\label{osim}
\ph(s\inv Ms)=s\inv\ph(M) s.
\eeq
\end{enumerate}

An nc-function $\ph$ on an nc-domain $D\subset \M^d$ is {\em analytic} if its restriction to $D\cap \M_n^d$ is analytic for every $n\in\mathbb{N}$.
\end{definition}
If $d=\infty$ the last statement should be interpreted to mean that $\ph$ is analytic with respect to the norm
$\|\cdot\|_n$ of Definition \ref{ncdomaininfty} on $D\cap \M_n^\infty$ for every $n$.

An nc-domain $D\subset \M^2$ is {\em symmetric} if $(M_2,M_1)\in D$ whenever $(M_1,M_2)\in D$.  
Clearly $\bb$ is a symmetric nc-domain.
If $\ph$ is an nc-function on a symmetric nc-domain $D\subset \M^2$, then $\ph$ is {\em symmetric} if $\ph(M_1,M_2)=\ph(M_2,M_1)$ for every $(M_1,M_2)\in D$.

\begin{definition} \label{ncanalonOmega}
If $D\subset \M^{d_1}$ and $ \Omega\subset \M^{d_2}$ are nc-domains, for $d_1,d_2\leq \infty$ then  an {\em nc-map} from $D$ to $\Omega$ is defined to be a map $F:D\to\Omega$ such that $F$ maps $D\cap  \M_n^{d_1}$ to  $\M_n^{d_2}$ for each $n\geq 1$ and $F$ respects direct sums and similarities, as in conditions \eqref{oadditive} and \eqref{osim}. 

If $\Omega$ is an nc-domain in $\M^\infty$ contained in a union of Banach spaces $\cup_{n\in\mathbb{N}} E_n$, $D$ is an nc-domain in $\M^d, \, d<\infty$, and $F: D\to \Omega$ is an nc-map 
then say that $F$ is {\em analytic} if, for each co-ordinate functional 
\[
f_j:\M^\infty\to \M^1: g=(g_0,g_1, \dots ) \mapsto g_j,
\]
the map $f_j\circ F$ is analytic for all $j\in\mathbb{N}$.
\end{definition}
An example of an nc-map from $\bb$ to an nc-domain $\Omega$ in $\M^\infty$ is the map $S$ described in Theorem \ref{realizn1}.

Operator-valued nc-functions will also be needed.  For Hilbert spaces $H$ and $K$
denote by $\mathcal{L}(H,K)$ the space of bounded linear operators from $H$ to $K$
with operator norm.  $\mathcal{L}(H,H)$ will be abbreviated to $\mathcal{L}(H)$.
An {\em $\mathcal{L}(H,K)$-valued nc-function} on an nc-domain $D$ in $\M^d$ 
is a function $\ph$ on $D$ such that
\begin{enumerate}
\item for $n\in\mathbb{N}$ and $M\in D\cap \M_m^d$, 
\[
\ph(M) \in \mathcal{L}(\C^n\otimes H,\C^n\otimes K);
\]
\item  for all $m,n\in\mathbb{N}$ and all $M\in D\cap \M_m^d$ and $N\in D \cap \M_n^d$,
\[
\ph(M\oplus N) = \ph(M) \oplus \ph(N)
\]
 modulo the natural identification of $(\C^m\otimes \mathcal{H}) \oplus (\C^n\otimes \mathcal{H})$ with $\C^{m+n} \otimes \mathcal{H}$ for any Hilbert space $\mathcal{H}$, and
\item  for any $m\in\mathbb{N}$, any $M\in D\cap \M_m^d$ and any invertible matrix $s\in\M_m$ such that $s\inv Ms \in D$,
\[
\ph(s\inv Ms) = (s\inv \otimes \1_K)\ph(M) (s\otimes \1_H).
\]
\end{enumerate}
The Hilbert space $H$ can be identified with $\mathcal{L}(\C,H)$ in the obvious way, so that
we may speak of $H$-valued nc-functions.

We shall denote the identity operator on any Hilbert space by $\1$.  Where
it is deemed particularly helpful to indicate the space we shall use subscripts;
thus $\1_n, \ \1_{\ell^2}$ are the identity operators on $\C^n, \ \ell^2$ respectively.

\section{Lurking isometries}\label{lurking}
A simple but powerful method in realization and interpolation theory is the use of {\em lurking isometries}:  if the gramians of two collections of vectors in Hilbert spaces are equal then there is an isometry that maps one collection to the other.  There is an nc version of the lurking isometry argument due to Agler and McCarthy; it is contained in the proof of \cite[Theorem 7.1]{JJ}.  

For an $\mathcal{L}(H,K)$-valued nc-function $f$ on an nc-domain $D$ in $\M^d$ (where $H,K$ are Hilbert spaces and $d\leq \infty$) define the {\em redundant subspace of $K$ for $f$}, denoted by $\Red (f)$, to be
\beq\label{codimf}
  \{\ga\in K : \C^n\otimes \ga \perp \bigvee_{x\in D\cap \M^d_n} \ran f(x) \mbox{ for all } n\in\mathbb{N} \}. 
\eeq
\begin{lemma}\label{nclurking}
Let $H, K_1$ and $K_2$ be Hilbert spaces and let $D$ be an nc-domain in $\M^d$ for some $d\leq\infty$.
Let $f$ be an $\mathcal{L}(H,K_1)$-valued nc-function and $g$ be an $\mathcal{L}(H,K_2)$-valued nc-function on $D$ such that, for all $n\geq 1$ and $x,y \in D\cap \M^d_n$,
\beq\label{f*f}
f(y)^*f(x) = g(y)^*g(x) \in \mathcal {L}(\C^n\otimes H).
\eeq
There exists a partial isometry $J:K_1\to K_2$ such that, for every positive integer $n$ and $x\in D\cap \M^d_n$,
\[
(\1_n\otimes J) f(x) = g(x).
\]
Moreover, if the dimensions of the redundant subspaces of $K_1$ and $K_2$ for $f$ and $g$ respectively
are equal then $J$ may be taken to be a unitary operator from $K_1$ to $K_2$.
\end{lemma}
\begin{proof}
Consider $x,y\in D\cap\M_n^d$ and an invertible $s\in\M_n$ such that $s\inv xs\in D$.  On replacing $x$ by $s\inv xs$ in equation \eqref{f*f} and invoking the fact that $f, g$ are nc-maps we have
\begin{align*}
f(y)^*(s\inv\otimes \1_{K_1}) f(x)(s\otimes \1_H) &=g(y)^*(s\inv\otimes \1_{K_2}) g(x)(s\otimes \1_H).
\end{align*}
Since the invertible matrices $s\inv$ with $\|s\|\, \|s\inv\|$ close to $1$ span all of $\M_n$ it follows that
\begin{align}\label{7.4}
f(y)^*(T\otimes \1_{K_1}) f(x)& = g(y)^*(T\otimes \1_{K_2}) g(x) \quad \in \mathcal{L}(\C^n\otimes H)
\end{align}
for all $T\in\M_n$.  Let $e_1,\dots e_n$ be the standard basis of $\C^n$ and apply equation \eqref{7.4} with $T=e_\ell e_k^*, \, k,\ell=1,\dots,n$, to deduce that, for any $\xi,\eta \in\C^n\otimes H$,
\beq\label{7.5}
\ip{(e_k^*\otimes \1_{K_1})f(x)\xi}{(e_\ell^*\otimes \1_{K_1})f(y)\eta}_{K_1} = \ip{(e_k^*\otimes \1_{K_2})g(x)\xi}{(e_\ell^*\otimes \1_{K_2})g(y)\eta}_{K_2}.
\eeq
Let
\begin{align*}
p_{k\xi x} &= (e_k^*\otimes \1_{K_1})f(x)\xi  \in K_1, \\
q_{k\xi x} &= (e_k^*\otimes \1_{K_2})g(x)\xi \in K_2
\end{align*}
and
\begin{align*}
\mathcal{P}_n &= \spa \{p_{k\xi x} : k\leq n, \xi\in\C^n\otimes H, \,  x\in D\cap \M^d_n\}\subset K_1, \\
\mathcal{Q}_n &= \spa \{q_{k\xi x} : k\leq n, \xi\in\C^n\otimes H,\,  x\in D\cap \M^d_n\}\subset K_2.
\end{align*}
Equation \eqref{7.5} states that
\[
\ip{p_{k\xi x}}{p_{\ell\eta y}}_{K_1} = \ip{q_{k\xi x}}{q_{\ell\eta y}}_{K_2}.
\]
It follows that there exists an isometry $L_n: \mathcal{P}_n \to \mathcal{Q}_n$ such that
\[
L_n p_{k\xi x} = q_{k\xi x}
\]
for all $k\leq n, \xi\in\C^n\otimes H$ and $x\in D\cap \M^d_n$.

We claim that both $(\mathcal{P}_n)$ and $(\mathcal{Q}_n)$ are increasing sequences of spaces, and
$L_m | \mathcal{P}_n = L_n$ when $n\leq m$.  Consider positive integers $n\leq m$ and regard $\C^n$ as  the span of the first $n$ standard basis vectors of $\ell^2$.  Let $k\leq n, \, \xi\in\C^n\otimes H$ and $x\in D\cap \M^d_n$.  For any choice of $x_0\in D\cap\M_{m-n}^d$ and $\xi_0\in\C^{m-n}\otimes H$,
\begin{align*}
p_{k(\xi\oplus\xi_0)(x\oplus x_0)} &= (e_k^*\otimes\1)f(x\oplus x_0)(\xi\oplus\xi_0)  \\
	&=(e_k^*\otimes\1)(f(x)\oplus f(x_0))(\xi\oplus\xi_0)  \\
	&=(e_k^*\otimes\1)\left(f(x) \xi \oplus f(x_0) \xi_0\right) \\
	&=(e_k^*\otimes\1)f(x) \xi \\
	&= p_{k\xi x}.
\end{align*}
Similarly $q_{k(\xi\oplus \xi_0)( x\oplus x_0)}= q_{k\xi x}$.
Hence $\mathcal{P}_n \subset \mathcal{P}_m$ and    $\mathcal{Q}_n \subset \mathcal{Q}_m$, while, for $k\leq n$,
\[
L_m p_{k\xi x} = q_{k\xi x} = L_n p_{k\xi x},
\]
 so that $L_m$ and $L_n$ agree on $\mathcal{P}_n $.

Let $\mathcal{P}, \mathcal{Q}$ be the closures in $K_1,\, K_2$ of $\cup_n \mathcal{P}_n,  \cup_n \mathcal{Q}_n$ respectively.  The isometries $L_n$ extend to an isometry $L:\mathcal{P} \to\mathcal{Q}$.  Extend $L$ further to a partial isometry $J:K_1\to K_2$.  Note that
\begin{align*}
K_1\ominus \mathcal{P} &= \{\ga\in K_1: \ip{(\eta^*\otimes \1_{K_1})f(x)\xi}{\ga}=0 \mbox{ for all } n\in\mathbb{N}, \xi,\eta \in \C^n, x\in D\cap \M^d_n \} \\
	&= \{\ga\in K_1: \ip{f(x)\xi}{(\eta\otimes \ga}=0 \mbox{ for all } n\in\mathbb{N}, \xi,\eta \in \C^n, x\in D\cap \M^d_n \} \\
	&= \{\ga\in K_1: \C^n\otimes \ga \perp \bigvee_{x\in D\cap \M^d_n} \ran f(x) \mbox{ for all } n\in\mathbb{N}\},
\end{align*}
which is the redundant subspace of $K_1$ for $f$.  Likewise
 $K_2\ominus \mathcal{Q}$ is the redundant subsace of $K_2$ for $g$.  Hence, if the dimensions of the two redundant subspaces spaces  are equal then the codimensions of $\mathcal{P}$ and $\mathcal{Q}$ in $K_1$ and $K_2$ respectively are equal, and consequently we can choose the partial isometry $J$ to be a unitary operator.
Whether or not the redundant subspaces have equal dimensions,  for any $n\in\mathbb{N}$ and for $x\in D\cap \M^d_n, \, \xi\in \C^n\otimes H$,
\begin{align*}
(\1_n\otimes J) f(x)\xi &= (\1_n\otimes J) \bigoplus_{k=1}^n(e_k^*\otimes \1_{K_1})f(x)\xi \\
	&=\bigoplus_{k=1}^n J(e_k^*\otimes \1_{K_1})f(x)\xi \\
	&=\bigoplus_{k=1}^n (e_k^*\otimes \1_{K_2})g(x)\xi \\
	&= q(x)\xi.
\end{align*}
Hence $(\1_n\otimes J) f(x)= g(x)$.
\end{proof}
Here is a simple property of nc-functions. 
\begin{proposition}\label{compose}
Let $H,K$ and $L$ be Hilbert spaces and let $D$ be an nc-domain in $\M^d$ for some $d\leq \infty$.
Let $f$ be an $\mathcal{L}(H,K)$-valued nc-function and let $g$ be an $\mathcal{L}(K,L)$-valued nc-function on $D$.  Then the function $gf$ defined by $(gf)(x)=g(x)f(x)$ for all $x\in D$ is an $\mathcal{L}(H,L)$-valued nc-function on $D$ and $\Red (g) \subset \Red(gf)$.  If $f(x)$ is an invertible operator for every $x\in D$ then $f(\cdot)\inv$ is an $\mathcal{L}(K,H)$-valued nc-function.
\end{proposition}
\begin{proof}
It is routine to show that $gf$ is an nc-function.  Suppose $\ga\in\Red(g)$: then for any $n\in\mathbb{N}, \, \xi\in\C^n\otimes K, \ \eta\in\C^n$ and $x\in D\cap \M^d_n$,
\[
\ip{\eta \otimes\ga}{g(x)\xi}_{\C^n\otimes L} =0.
\]
In particular this holds when $\xi=f(x)\xi'$ for any $\xi' \in \C^n\otimes H$, which implies that $\ga\in\Red(gf)$.
\end{proof}
\begin{proposition}\label{boldAofD}
The nc-disc algebra $\mathbf{A}(\D)$ is a union of Banach spaces contained in $\M^\infty$ with respect to the norms
\[
\|g\| = \sup_{z\in\D} \|g(z)\|  \quad \mbox{ for all }g\in \mathbf{A}(\D) \cap \M_n^\infty  \mbox{ and all } n\in\mathbb{N}
\]
when the function $g\in \mathbf{A}(\D)$ is identified with its sequence of Taylor coefficients.
\end{proposition}
\begin{proof}
The space
\beq\label{defAn}
A_n(\D) \df \mathbf{A}(\D) \cap \M_n^\infty,
\eeq
the {\em $n\times n$-disc algebra}, is clearly a Banach space for the supremum norm, and this norm 
induces a stronger topology than the topology of pointwise convergence of sequences of Taylor coefficients,
which is the product topology on $\M^\infty$ restricted to $A_n(\D)$.  The supremum norm is also 
unitarily invariant.  Hence $\mathbf{A}(\D)$ is a union of Banach spaces contained in $\M^\infty$
in the sense of Definition \ref{ncdomaininfty}.
\end{proof}

$\mathbf{A}(\D)$ has the structure of an {\em operator space}, but we shall not use this fact.

\section{Linear fractional maps}\label{linfrac}
 For any block matrix
\beq\label{defP}
p = \bbm p_{11}  & p_{12} \\ p_{21} &p_{22}  \ebm
\eeq
we shall denote by $\F_p^\ell$ the {\em lower linear fractional transformation}
\beq\label{defLFT}
\F_p^\ell (X) = p_{22} + p_{21} X (\1-  p_{11} X)^{-1}p_{12}
\eeq
whenever the formula is meaningful.  
For example, when $p_{ij}$ is an $m_i\times n_j$ matrix, it is defined for every $n_1\times m_1$
matrix $X$ such that $1 - p_{11}X$ is invertible, and then
$\F_p^\ell (X)$ is an $m_2\times n_2$ matrix.  More generally, if $\mathcal{H}_i, \mathcal{K}_i$ are Hilbert spaces for $i=1,2$ and $p$ is a block operator matrix from $\mathcal{K}_1\oplus\mathcal{H}_2$ to $\mathcal{H}_1\oplus     \mathcal{K}_2$ and $X$ is a bounded operator from $\mathcal{H}_1$ to $\mathcal{K}_1$ such that $1 - p_{11}X$ is invertible on $\mathcal{H}_1$ then $\F_p^\ell  (X)$  is defined and is an operator from $\mathcal{H}_2$ to $\mathcal{K}_2$.

We shall also define the {\em upper linear fractional transformation}
\beq\label{defFu}
\F_p^u (X) = p_{11}+p_{12}X(1-p_{22}X)\inv p_{21}.
\eeq
 \black

The following results are standard.
\begin{lemma}\label{LFbasic}
For any matrices or operators $p, X$ such that $\F_p^\ell  (X) $
is defined
\begin{eqnarray}\label{realization}
\1 - \F_p^\ell  (X)^* \F_p^\ell  (X) 
&=&   p_{12}^* \left(\1 - X^* p^*_{11}\right)^{-1}
(\1- X^* X) \left(\1 - p_{11} X \right)^{-1} p_{12} \nonumber\\
&~& \hspace{-3cm}
+\bbm p_{12}^* \left(\1-Xp^*_{11}\right)^{-1} X^*  & \1 \ebm  
(\1- p^* p) 
\bbm X \left(\1 - p_{11} X \right)^{-1} p_{12}  \\ \1 \ebm.
\end{eqnarray} 
Furthermore,  if $p, X$ are contractions then
\beq\label{normfpx}
\|\F_p^\ell  (X) \| \leq \|p\|.
\eeq
\end{lemma}
Of course analogous results hold for $\mathcal{F}_{ p}^u$.
\begin{proof}
The identity \eqref{realization} may be verified by straightforward expansion.
Since $\1-X^*X \geq 0$ the identity implies that
\[
\1 - \F_p^\ell (X)^* \F_p^\ell  (X) \geq   \bbm p_{12}^* \left(\1-Xp^*_{11}\right)^{-1} X^*  & \1 \ebm  
(\1- p^* p) 
\bbm X \left(\1 - p_{11} X \right)^{-1} p_{12}  \\ \1 \ebm.
\]
Hence, since $\1-p^*p \geq (1-\|p\|^2)\1$, for any vector $\xi$,
\begin{align*}
\ip{(\1 - \F_p^\ell  (X)^* \F_p^\ell  (X))\xi}{\xi} &\geq \ip{(1-\|p\|^2)\bbm X \left(\1 - p_{11} X \right)^{-1} p_{12}  \\ \1 \ebm\xi}{\bbm X \left(\1 - p_{11} X \right)^{-1} p_{12}  \\ \1 \ebm\xi} \\
	&\geq (1-\|p\|^2)\|\xi\|^2.
\end{align*}
The inequality \eqref{normfpx} follows.
\end{proof}

There are also {\em  graded linear fractional maps}, which map $n\times n$ matrices to $n\times n$ matrices.
\begin{definition} \label{linfracnc}
Let $\mathcal{H}_i, \mathcal{K}_i$ be Hilbert spaces for $i=1,2$ and let $p$ be a block operator matrix from $\mathcal{K}_1\oplus\mathcal{H}_2$ to $\mathcal{H}_1\oplus     \mathcal{K}_2$.  The {\em graded lower linear fractional map with matrix $p$} is defined to be the map $\mathcal{F}_{\1\otimes p}^\ell$ which, for $n\geq 1$, maps $X\in\mathcal{L}(\C^n\otimes \mathcal{H}_1, \C^n\otimes\mathcal{K}_1)$ such that $\1-(\1_n\otimes p_{11})X$ is invertible to
\[
\mathcal{F}_{\1\otimes p}^\ell(X)= \1_n\otimes p_{22} + (\1_n\otimes p_{21})X\left(\1-(\1_n\otimes p_{11})X\right)\inv (\1_n\otimes p_{12}).
\]
Similarly we define the {\em graded upper linear fractional map}
\beq\label{upLFnc}
\mathcal{F}_{\1\otimes p}^u (X) = \1_n\otimes p_{11} + (\1_n\otimes p_{12})X\left(\1-(\1_n\otimes p_{22})X\right)\inv (\1_n\otimes p_{21}).
\eeq
\end{definition}


Observe that $\mathcal{F}_{\1\otimes p}^\ell(X)$ is an operator from $\C^n\otimes\mathcal{H}_2$ to $\C^n\otimes\mathcal{K}_2$ for each $n$.  Likewise $\mathcal{F}_{\1\otimes p}^u(X)$ is defined for suitable operators $X:\C^n\otimes\mathcal{K}_2\to \C^n\otimes\mathcal{H}_2$ and is an operator from $\C^n\otimes\mathcal{K}_1$ to $\C^n\otimes\mathcal{H}_1$ for each $n$.

The function $\mathcal{F}_{\1\otimes p}^u$ enjoys some properties of nc type.  Its domain is the set
\[
D= \bigcup_{n=1}^\infty \{X\in \mathcal{L} (\C^n\otimes \mathcal{K}_2, \C^n\otimes \mathcal{H}_2) : \1-(\1_n\otimes p_{22})X \mbox{ is invertible} \}.
\]
\begin{proposition}\label{partnc}
Let $p$ be the block operator matrix from $\mathcal{K}_1\oplus\mathcal{H}_2$ to $\mathcal{H}_1\oplus     \mathcal{K}_2$ given by equation \eqref{defP}.  Its domain
$D$ is closed under direct sums, and, for $X,Y \in D$,
\[
\mathcal{F}_{\1\otimes p}^u(X\oplus Y) = \mathcal{F}_{\1\otimes p}^u(X) \oplus \mathcal{F}_{\1\otimes p}^u(Y).
\]
Moreover, if $X\in D\cap \mathcal{L} (\C^n\otimes \mathcal{K}_2, \C^n\otimes \mathcal{H}_2)$ and $s\in\M_n$ is an invertible matrix then $(s\inv\otimes \1_{\mathcal{H}_2}) X (s\otimes \1_{\mathcal{K}_2}) \in D$ and
\[
\mathcal{F}_{\1\otimes p}^u((s\inv\otimes \1_{\mathcal{H}_2}) X (s\otimes \1_{\mathcal{K}_2})) =
(s\inv\otimes \1_{\mathcal{H}_1}) \mathcal{F}_{\1\otimes p}^u( X)(s\otimes \1_{\mathcal{K}_1}).
\]
\end{proposition}
The proof is straightforward.

\section{A realization theorem}\label{realize}
In this section we show that every bounded symmetric nc-function on the biball factors through a certain nc-domain $\Omega$ in $\mathcal{M}^\infty$
and is thereby expressible by a linear fractional realization formula.

$\M^\infty$ is natually identified with the space $\M[[z]] =\cup_n\M_n[[z]]$ of formal power series over $\M$ in the indeterminate $z$.  For $n\geq 1$ the element $g=(g^0,g^1,\dots) \in \M_n^\infty$ corresponds to the series $\sum_{j\geq 0} g^j z^j \in \M_n[[z]]$. With this understanding the `functional calculus map' $\Theta_T$ on (a subset of) $\M^\infty$ corresponding to an operator $T$ on a Hilbert space $H$ is given by
\beq\label{defTheta}
\Theta_T(g) = \sum_{j=0}^\infty g^j \otimes T^j \in \mathcal{L}(\C^n\otimes H),
\eeq
 whenever the series converges in an appropriate sense.  In the present context it is enough that the series
in equation \eqref{defTheta}  converge in the sense of C\'esaro summability of the partial sums of the series 
in the operator norm.
As is customary, $\Theta_T(g)$ will also be denoted by $g(T)$ when it exists.

\begin{theorem}\label{realizn}
  There exists an nc-domain $\Omega$ in $\M^\infty$ such that the map $S:\bb\to \M^\infty$ defined by
\beq\label{defS}
S( x)= (u,v^2,vuv,vu^2v,\dots),
\eeq
where
\beq\label{defuv}
u=\frac{x^1+x^2}{2}, \qquad v=\frac {x^1-x^2}{2},
\eeq
has the following three properties.
\begin{enumerate}
\item $S$  is an analytic nc-map from $\bb$ to $\Omega$;
\item  for every $g\in\Omega$ and every contraction $T$ the operator $g(T)$ exists and $\|g(T)\|<1$;
 \item for every symmetric nc-function $\ph$ on the biball bounded by $1$ in norm there exists an analytic  nc-function $\Phi$ on $\Omega$ such that $\|\Phi(g)\|\leq 1$ for every $g\in\Omega$ and $\ph=\Phi\circ S$.
\end{enumerate}

Moreover $\Phi$ can be realized as follows.  There exist a unitary operator $U$ on $\ell^2$ and a contractive operator
\[
p=\bbm p_{11}& p_{12}\\p_{21}&p_{22} \ebm:\C\oplus \ell^2 \to \C\oplus \ell^2
\]
such that, for $n\geq 1$ and $g\in\Omega\cap \M_n^\infty$,
\begin{align}\label{defPhi}
\Phi(g)&= \mathcal{F}_{\1\otimes p}^u(g(U)) \\
	&= p_{11}\1_n + (\1_n\otimes p_{12})g(U)\left(\1-(\1_n\otimes p_{22})g(U)\right)\inv (\1_n\otimes p_{21}). \notag
\end{align}
\end{theorem}

\begin{proof}  The existence of models of bounded nc-functions on the polyball is a major result of \cite{JJ}; we shall combine this result with a symmetrization argument.

Let $\Omega$ be the open unit ball of the nc-disc algebra $\mathbf{A}(\D)$ of Example \ref{exAofD}.  More precisely, $\Omega$ is the union of the open unit balls of the Banach spaces $A_n(\D)= \mathbf{A}(\D) \cap \M_n^\infty$ for $n\geq 1$.  By Proposition \ref{boldAofD} $\mathbf{A}(\D)$ is a union of Banach spaces contained in $\M^\infty$, and it is easy to see that $\Omega$ is an nc-domain in $\M^\infty$.

To prove (2) consider any $g\in\Omega\cap\M_n^\infty$ and any contractive operator $T$ on a Hilbert space $\mathcal{H}$.  For $k\geq 0$ let $h_k(z)$ be the arithmetic mean of the $k+1$ Taylor polynomials 
\[
 g^0+g^1 z+\dots+g^r z^r, \qquad r=0,1,\dots, k
\]
of $g\in A_n(\D)$.
By Fej\'er's theorem $h_k$ converges uniformly on $\D^-$ to $g$.  By von Neumann's inequality $ (h_k(T))_{k\geq 1}$ is a Cauchy sequence with respect to the operator norm, and so $g(T)=\Theta_T(g)$ is defined to be the limit of the sequence $(h_k(T))$ in $\mathcal{L}(\C^n\otimes\mathcal{H})$.  Since $\|h_k\|_\infty \to \|g\|_\infty <1$, it follows that $\|g(T)\|<1$.

For (1) consider $x\in\bb\cap \M_n^\infty$: we must prove that $S(x)\in\Omega$. If $S(x)$ is identified with its generating function $S(x)(z)$ then, since $\|v\|<1$,
\begin{align*}
S(x)(z) &= u+v^2 z+ vuv z^2+vu^2v z^3+\dots \\
	&= u+ vz(\1_n-uz)\inv v \\
	&= \frac{ x^1 +  x^2}{2} +
	 \frac{ x^1 -  x^2}{2}z  
	\left(\1_n -\frac{x^1+x^2}{2}z\right)\inv \frac{ x^1 -  x^2}{2} 
\end{align*}
Clearly $S(x) \in \mathbf{A}(\D)$.
Let
\[
Q(x) =\bbm u&v\\v&u\ebm.
\]
Then
\[
Q(x) = \frac{1}{2}\bbm  \ds x^1+x^2 & \ds x^1-x^2\\ x^1-x^2 &x^1+x^2 \ebm
	=\frac{1}{\sqrt{2}}\bbm \1 & \1\\ \1 & -\1 \ebm \bbm x^1&0\\0&x^2 \ebm \frac{1}{\sqrt{2}}\bbm \1 & \1\\ \1 & -\1 \ebm
\]
and hence
\[
\|Q(x)\|=\max\{\|x^1\|,\|x^2\|\} =\|x\| < 1.
\]
Since 
\[
S(x)(z) = \F_{Q(x)}(z\1_n)
\]
it follows from Lemma \ref{LFbasic}  that 
\[
\|S(x)(\cdot)\|_{A_n(\D)} \  \leq \ \|x\| \ < \ 1,
\]
 and so $S(x)\in\Omega$.

If $x\in\bb\cap \M_n^2$ then $S(x) \in\Omega \cap \M_n^\infty$ for each $n\geq 1$.  Moreover $S$ respects direct sums and similarities: if $x\in\bb\cap \M_m^2$ and $y\in\bb\cap \M_n^2$ then, for $z\in\D$,
\begin{align*}
S(x\oplus y)(z) &= \frac{x^1\oplus y^1 + x^2\oplus y^2}{2}  +\\
	&\hspace*{-0.5cm} \frac{x^1\oplus y^1 - x^2\oplus y^2}{2} z\left( \1_{m+n}-\frac{x^1\oplus y^1 + x^2\oplus y^2}{2} z\right)\inv \frac{x^1\oplus y^1 - x^2\oplus y^2}{2} \\
	&= S(x)(z)\oplus S(y)(z),
\end{align*} 
while if $s$ is an invertible matrix in $\M_m$ such that $s\inv xs\in\bb$ then, for $z\in\D$,
\begin{align*}
S(s\inv xs)(z)&= \frac{s\inv x^1s + s\inv x^2s}{2} +\\
	& \hspace* {-1cm} \frac{s\inv x^1s - s\inv x^2s}{2}z  
	\left(\1_m -\frac{s\inv x^1s+ s\inv x^2s}{2}z\right)\inv \frac{s\inv x^1s - s\inv x^2s}{2} \\
	&= s\inv S(x)(z) s.
\end{align*}
Hence $S$ is an nc-map from $\bb$ to $\Omega$.  It is analytic since the restriction of $S$ mapping $\bb\cap \M_n^2$ to $\Omega \cap \M_n^\infty\subset A_n(\D)$ is an analytic Banach-space-valued map for each $n$.  We have proved (1).

Let $\ph$ be a symmetric nc-function on $\bb$ bounded by $1$ in norm.  By \cite[Theorem 6.5]{JJ} $\ph$ has an nc-model; that is, there is a pair $(P,\mm )$ where $P=(P^1,P^2)$ is an orthogonal decomposition of $\ell^2$ (so that $P^1+P^2=\1_{\ell^2}$), $\mm $ is an $\ell^2$-valued nc-function on $\bb$ and
\beq\label{modeleq}
\1_n-\ph(y)^*\ph(x) = \mm (y)^*(\1_{\C^n\oplus \ell^2}-y_P^*x_P)\mm (x)
\eeq
for all $x,y\in\bb\cap \M_n^2$.  Here $x_P$ denotes $x^1\otimes P^1+x^2\otimes P^2$, an operator on $\C^n\otimes \ell^2$.  

Since
\begin{align*}
\1-y_P^*x_P &= \1- (\sum_j y^j\otimes P^j)^*(\sum_ix^i\otimes P^i) \\
	&= \1-\sum_i (y^{i*}x^i \otimes P^i) \\
	&= \sum_i(\1-y^{i*}x^i)\otimes P^i,
\end{align*}
 equation \eqref{modeleq} can also be written (in the case that $x,y\in \bb\cap\M^2_n$)
\begin{align}\label{newmodeleq}
\1_n -\ph(y)^*\ph(x) &=\mm(y)^* \sum_i((\1_n-y^{i*}x^i)\otimes P^i) \mm(x)\notag\\
	&= \sum_i\mm(y)^*(\1_n\otimes P^i)\left((\1_n-y^{i*}x^i)\otimes \1_{\ell^2}\right) (\1_n\otimes P^i)\mm(x)\notag\\
	&= \sum_{i=1}^2 \mm^i(y)^* ((\1_n-y^{i*}x^i)\otimes\1_{P^i\ell^2})\mm^i(x)
\end{align}
where, for $i=1,2$,
\[
\mm^i(x)\df (\1_n\otimes P^i)\mm(x) \in \mathcal{L}(\C^n, \C^n\otimes P^i\ell^2).
\]

Let $H_i=P^i\ell^2$ for $i=1,2$.
We claim that $\mm^i$ is an  $H_i$-valued nc-function on $\bb$.  Certainly $\mm^i(x)\in\mathcal{L}(\C^n,\C^n\otimes H_i)$ for $x\in \bb\cap \M_n^2$.  If $x, y \in\bb$ are $n$-square and $m$-square respectively then
\[
\mm^i(x\oplus y)=(\1_{n+m} \otimes P^i)\mm(x\oplus y) = ((\1_{n}\oplus \1_m) \otimes P^i)(\mm(x)\oplus \mm( y)) = \mm^i(x) \oplus \mm^i(y).
\]
Furthermore, if $s\in\M_n$ is invertible and $s\inv x s$ belongs to $\bb$ for some $x\in \bb\cap\M^2_n$, then
\begin{align*}
\mm^i(s\inv xs)&=(\1_n\otimes P^i)\mm(s\inv xs) \\
	&=(\1_n\otimes P^i)(s\inv \otimes \1_{\ell^2})\mm(x)s \\
	&=(s\inv \otimes \1_{H_i}) (\1_n\otimes P^i) \mm(x)s \\
	&=(s\inv \otimes \1_{H_i})  \mm^i(x)s.
\end{align*}
Thus $\mm^i$ is an $H_i$-valued nc-function on $\bb$ as claimed.

Since $\ph$ is symmetric we may interchange $y^1$ and $y^2$, $x^1$ and $x^2$ in equation \eqref{newmodeleq} to obtain
\beq\label{6.2}
\1_n-\ph(y)^*\ph(x)= \tilde\mm^1(y)^*\left((\1_n-y^{2*}x^2)\otimes \1_{H_1}\right) \tilde\mm^1(x) + \tilde \mm^2(y)^*\left( (\1_n-y^{1*}x^1)\otimes \1_{H_2}\right) \tilde \mm^2(x)
\eeq
where, for any function $\psi$ on $\bb$, $\tilde\psi (x^1,x^2)$ denotes $\psi(x^2,x^1)$.  Notice that $\tilde\mm^i$ is also an $H_i$-valued nc-function on $\bb$.

Average equations \eqref{newmodeleq}, \eqref{6.2} to deduce that
\beq\label{6.3}
\1_n-\ph(y)^*\ph(x)= w(y)^*\left((1-y^{1*}x^1)\otimes \1_{\ell^2}\right) w(x) +  \tilde w(y)^*\left( (1-y^{2*}x^2)\otimes \1_{\ell^2}\right) \tilde w(x)
\eeq
 for $x,y \in \bb$, where
\beq\label{defw}
w(x) = \frac{1}{\sqrt{2}} \bbm \mm^1(x) \\ \tilde \mm^2(x)\ebm :\C^n\to(\C^n\otimes H_1)\oplus (\C^n\otimes H_2)= \C^n\otimes\ell^2.
\eeq
and so
\[
\tilde w(x) = \frac{1}{\sqrt{2}} \bbm \tilde \mm^1(x) \\  \mm^2(x) \ebm:\C^n\to\C^n\otimes\ell^2.
\]
Since $\mm^i, \tilde \mm^i$ are $H_i$-valued nc-functions on $\bb$, the functions $w$ and $\tilde w$ are  $\ell^2$-valued nc-functions.

In equation \eqref{6.3} interchange $x^1,x^2$ (but not $y^1,y^2$) and use the symmetry of $\ph$ to deduce that
\begin{align*}
w(y)^*\left( (\1_n-y^{1*}x^1) \right. &\left. \otimes \1_{\ell^2}\right) w(x) +\tilde w(y)^* \left((\1_n-y^{2*}x^2)\otimes \1_{\ell^2}\right) \tilde w(x) = \\
	&  w(y)^*\left((\1_n-y^{1*}x^2) \otimes \1_{\ell^2}\right) \tilde w(x) +  \tilde w(y)^* \left((\1_n-y^{2*}x^1)\otimes \1_{\ell^2}\right)  w(x).
\end{align*}
Rearrangement of this equation yields
\begin{align*}
   w(y)^*w(x) &+ \tilde w(y)^* \tilde w(x) - w(y)^* \tilde w(x)- \tilde w(y)^*w(x)  =\\
	& w(y)^*( y^{1*}x^1\otimes \1_{\ell^2})w(x)+ \tilde w(y)^* ( y^{2*}x^2 \otimes \1_{\ell^2})\tilde w(x)  \\
	&\hspace*{1cm}   - w(y)^* ( y^{1*}x^2\otimes \1_{\ell^2})\tilde w(x)- \tilde w(y)^*( y^{2*}x^1 \otimes \1_{\ell^2})w(x).
\end{align*}
Both sides of the equation factor:
\begin{align} \label{7.2}
 \left(w(y)^* - \tilde w(y)^*  \right)  \left(w(x) - \tilde w(x) \right) &= 
	\quad\left( w(y)^*( y^{1*}\otimes\1_{\ell^2}) - \tilde w(y)^*( y^{2*}\otimes \1_{\ell^2}) \right)\notag \\
	&\times \left( ( x^1\otimes \1_{\ell^2}) w(x) - ( x^2\otimes \1_{\ell^2})\tilde w(x)   \right)
\end{align}

Since both $w$ and $\tilde w$ are $\ell^2$-valued nc-functions on $\bb$, so are $w-\tilde w$ and the function
\[
g(x)= (x^1\otimes \1_{\ell^2}) w(x) - ( x^2\otimes \1_{\ell^2})\tilde w(x).
\]

We can assume that the redundant spaces of both $w-\tilde w$ and $g$ are infinite-dimensional.  To see this replace the nc-model $(P,\mm)$ of $\ph$ by the model $(Q,\psi)$ with model space $\ell^2\oplus\ell^2$ (which may be identified with $\ell^2$) which is trivial on the first copy of $\ell^2$ and agrees with $(P,\mm)$ on the second copy.  More precisely, $Q^1=0\oplus P^1, \, Q^2=0\oplus P^2$ and 
\[
\psi(x) = 0 \oplus \mm(x): \C^n \to (\C^n\otimes \ell^2) \oplus (\C^n\otimes \ell^2)
\]
for $x\in \bb\cap\M^2_n$.  Then
\begin{align*}
\psi(y)^*(\1_{\C^n\otimes(\ell^2\oplus\ell^2)} -y_Q^*x_Q)\psi(x) &=\bbm 0&\mm(y)^*\ebm
	\left( \1 -\diag(0, y_P^*x_P)\right) \bbm 0\\ \mm(x) \ebm \\
	&= \mm(y)^*(\1-y_P^*x_P)\mm(x) \\
	&= \1_n-\ph(y)^*\ph(x)
\end{align*}
and so $(Q,\psi)$ is a model of $\ph$.  It is easy to see that $\psi$ is an $\ell^2\oplus\ell^2$-valued nc-function on $\bb$.  Now if $w^\sharp$ is the analog of $w$ defined with $(Q,\psi)$ instead of $(P,\mm)$ then for $x\in\bb\cap\M^2_n$,
\[
w^\sharp (x) = 0 \oplus w(x): \C^n \to (\C^n \otimes \ell^2)\oplus (\C^n \otimes \ell^2)
\]
and the redundant space of $\ell^2\oplus\ell^2$ for $w^\sharp- \tilde w^\sharp$  contains $\ell^2 \oplus \{0\}$: for $\xi,\eta\in\C^n$ and $x\in\bb\cap\M^2_n$ and $\zeta\in\ell^2$,
\begin{align*}
\ip{\eta\otimes(\zeta\oplus 0)}{(w^\sharp(x)-\tilde w^\sharp(x))\xi}_{\C^n \otimes(\ell^2\oplus\ell^2)} &=
	\ip{(\eta\otimes\zeta)\oplus 0_{\C^n\otimes \ell^2}}{0_{\C^n\otimes \ell^2} \oplus (w(x)-\tilde w(x))\xi} \\
	&= 0
\end{align*}
and so $\C^n\otimes (\zeta\oplus 0) \perp \ran \psi(x)$ for all $n\in\mathbb{N}$ and $x\in\bb\cap\M^2_n$.  Similarly $\ell^2 \oplus \{0\}$ is contained in the redundant subspaces of $\ell^2\oplus\ell^2$ for $w^\sharp$ and for $g^\sharp$.

With the assumption of infinite-dimensional redundant subspaces of $w$ and $g$, by Lemma \ref{nclurking} 
there exists a unitary operator $U$ on $\ell^2$ such that for all $n \geq 1$ and all $x\in \bb \cap \M_n$,
\beq\label{8.1}
w(x) -\tilde w(x)= ( \1_n\otimes U)\left( ( x^1\otimes \1_{\ell^2}) w(x) -( x^2\otimes \1_{\ell^2})\tilde w(x)\right).
\eeq
and hence
\beq \label{8.2}
\left(\1- ( x^1\otimes U) \right)w(x) = \left(\1-( x^2\otimes U)\right) \tilde w(x).
\eeq
We wish to rewrite the model relation \eqref{6.3} incorporating the equation \eqref{8.2}.  To make it more concise let us introduce the abbreviations
\begin{align*}
\vv(x)&= \left(\1- ( x^1\otimes U) \right)w(x): \C^n\to \C^n\otimes \ell^2, \\
X^j&= x^j\otimes U \in \mathcal{L}(\C^n\otimes \ell^2), \\
Y^j&= y^j\otimes U \in \mathcal{L}(\C^n\otimes \ell^2)
\end{align*}
for $j=1,2$.   It is straightforward to check that $\vv$ is an $\ell^2$-valued nc-function on $\bb$.
Equation \eqref{8.2} states that $\vv(x)$ is symmetric in $(x^1,x^2)$, and we have
\[
w(x) = (\1-X^1)\inv \vv(x), \quad \tilde w(x) =(\1-X^2)\inv \vv(x)
\]
and
\[
Y^{1*}X^1 =   y^{1*}x^1 \otimes \1_{\ell^2}
\]
In terms of $\vv$  and $X^j, \ Y^j$ the model relation \eqref{6.3} can be written
\begin{align*}
\1_n-\ph(y)^*\ph(x) &= 
	 \vv(y)^*(\1-Y^{1*})\inv(\1-Y^{1*}X^1)(\1-X^1)\inv \vv(x) \\  &\quad + \vv(y)^*(\1-Y^{2*})\inv (\1-Y^{2*}X^2)(\1-X^2)\inv \vv(x).
\end{align*}
 Now
\begin{align*}
(\1-Y^{1*})\inv (\1-Y^{1*}X^1)& (\1-X^1)\inv \\
	&= (\1-Y^{1*})\inv(\1-X^1)\inv - (\1-Y^{1*})\inv Y^{1*}X^1 (\1-X^1)\inv \\
	&=(\1-Y^{1*})\inv(\1-X^1)\inv - \left((\1-Y^{1*}) -\1\right)\left(  (\1-X^1)-\1\right)\inv \\
	&=(\1-Y^{1*})\inv + (\1-X^1)\inv-\1.
\end{align*}
Hence equation \eqref{6.3} becomes
\begin{align}\label{10.1}
\1_n-\ph(y)^*\ph(x) &= \vv(y)^*\left[(\1-Y^{1*})\inv + (\1-X^1)\inv-\1 +\right. \notag \\
	&\hspace*{1.5cm}\left .  (\1-Y^{2*})\inv + (\1-X^2)\inv-\1 \right]\vv(x)  \notag\\
	&=\vv(y)^*\left[A(x) + A(y)^* \right]\vv(x)
\end{align}
where
\begin{align}\label{10.2}
A(x) &= (\1-X^1)\inv+(\1-X^2)\inv-\1 \notag \\
	&=(\1-( x^1\otimes U))\inv + (\1-( x^2\otimes U))\inv -\1 \quad \in \mathcal{L}(\C^n \otimes \ell^2).
\end{align}
It is easy to verify that $A$ is an $\mathcal{L}(\ell^2)$-valued nc-function on $\bb$.
Since
\[
A(x)+A(y)^*=\half (\1+A(y))^*(\1+A(x)) - \half (\1-A(y))^*(\1-A(x)),
\]
equation \eqref{10.1} implies that, for any $x,y\in\bb$,
\[
\1_n-\ph(y)^*\ph(x) = \half \vv(y)^*(\1+A(y))^*(\1+A(x))\vv(x) -\half \vv(y)^*(\1-A(y))^*(\1-A(x)) \vv(x).
\]
The last equation can also be written
\[
\bbm \1_n \\ \frac{1}{\sqrt{2}}(\1-A(y))\vv(y) \ebm^*\bbm \1_n \\ \frac{1}{\sqrt{2}}(\1-A(x))\vv(x) \ebm =
\bbm \ph(y)\\ \frac{1}{\sqrt{2}}(\1+A(y))\vv(y) \ebm^* \bbm \ph(x)\\ \frac{1}{\sqrt{2}}(\1+A(x))\vv(x) \ebm.
\] 
Since both $\vv$ and $A$ are nc-functions, the maps
\[
x\in\bb \cap \M^2_n \mapsto \bbm \1_n \\ \frac{1}{\sqrt{2}}(\1\pm A(x))\vv(x) \ebm \in\mathcal{L}(\C^n,\C^n\otimes (\C\oplus \ell^2))
\]
are $(\C\oplus\ell^2)$-valued nc-functions.  Hence by Lemma \ref{nclurking} there exists a contraction
\beq\label{defT}
T\df \bbm a&B\\C&D \ebm : \C\oplus \ell^2 \to \C\oplus\ell^2
\eeq
such that, for $n\geq 1$ and $x\in \bb\cap \M^2_n$,
\beq\label{11.1}
\bbm \ph(x) \\ \frac{1}{\sqrt{2}} (\1+A(x))\vv(x) \ebm = \bbm a\1_n & \1_n \otimes B \\ \1_n \otimes C  & \1_n\otimes D \ebm  \bbm \1_n \\  \frac{1}{\sqrt{2}} (\1-A(x))\vv(x) \ebm.
\eeq

We need a simple matrix identity.
\begin{lemma}\label{simple}
Let $Z_1,Z_2 \in\M_n$ and suppse that $Z_1,Z_2$ and $Z_1+Z_2$ are all invertible.  Then $Z_1\inv+Z_2\inv$ is invertible and
\[
4(Z_1\inv+Z_2\inv)\inv = Z_1+Z_2 - (Z_1-Z_2)(Z_1+Z_2)\inv(Z_1-Z_2).
\]
\end{lemma}
\begin{proof}
\begin{align*}
(Z_1-Z_2)(Z_1+Z_2)\inv(Z_1-Z_2)&=(Z_1+Z_2-2Z_2)(Z_1+Z_2)\inv(Z_1-Z_2)\\
	&=Z_1-Z_2-2Z_2(Z_1+Z_2)\inv(Z_1-Z_2)\\
	&=Z_1-Z_2-2Z_2(Z_1+Z_2)\inv(2Z_1-(Z_1+Z_2))\\
	&=Z_1-Z_2-4Z_2(Z_1+Z_2)\inv Z_1 + 2Z_2\\
	&=Z_1+Z_2 - 4(Z_1\inv +Z_2\inv)\inv.
\end{align*}
\end{proof}
Resume the proof of Theorem \ref{realizn}.   From the definition \eqref{10.2} of $A(x)$ and Lemma \ref{simple} with $Z_j=\1-X^j$,
\begin{align}\label{11.3}
\frac{ \1-A(x)}{\1+A(x)} &= -\1+\frac{2}{\1+A(x)}\notag \\
	&=-\1+\frac{2}{(\1-X^1)\inv + (\1-X^2)\inv} \notag\\
	&=-\1+\half\{2 \1-X^1-X^2)-(X^1-X^2) (2\1 -X^1-X^2)\inv (X^1-X^2)\}\notag\\
	&= -\frac{X^1+X^2}{2} - \frac{X^1-X^2}{2}\left(\1-\frac{X^1+X^2}{2}\right)\inv\frac{X^1-X^2}{2}\notag \\
	&=- \left(\frac{x^1+x^2}{2} \otimes U+ \frac{x^1-x^2}{2}\otimes U\left(\1-\frac{x^1+x^2}{2}\otimes U\right)\inv\frac{x^1-x^2}{2}\otimes  U\right)
\end{align}
which is an operator on $ \C^n\otimes \ell^2$ when $x\in\M_n^2$.   

Recall the notations $u=\half(x^1+x^2), \, v=\half(x^1-x^2)$ and $S(x)=(u,v^2,vuv,vu^2v,\dots)\in\M^\infty$.  We have, for any $x\in\bb$,
\begin{align*}
\Theta_U(S(x)) &= u\otimes \1_{\ell^2}+ v^2\otimes U+vuv\otimes U^2+vu^2v\otimes U^3+\dots \\
	&=  u\otimes \1_{\ell^2}+ (v\otimes U ) (\1 - u\otimes U)\inv (v\otimes \1_{\ell^2}),
\end{align*}
so that equation \eqref{11.3} becomes
\beq\label{concise}
\frac{ \1-A(x)}{\1+A(x)} = - (\1_n\otimes U)\Theta_U(S(x)).
\eeq

Next combine equations \eqref{concise} and \eqref{11.1} to obtain a realization formula for $\ph$.
To this end write 
\[
\vv^\flat(x)= \frac{1}{\sqrt{2}}(\1+A(x))\vv(x) \in \mathcal{L}(\C^n, \C^n\otimes \ell^2),
\]
so that 
\begin{align*}
\frac{1}{\sqrt{2}}(\1-A(x))\vv(x) &= \frac{\1-A(x)}{\1+A(x)}\vv^\flat(x) \\ 
	&=   - (\1_n\otimes U)\Theta_U(S(x))\vv^\flat(x).
\end{align*}
Equation \eqref{11.1} can thus be written as the pair of relations
\begin{align}\label{pair1}
\ph(x) &= a\1_n +(\1_n \otimes B)\frac{1}{\sqrt{2}}(\1-A(x))\vv(x) \notag\\
	&=a\1_n +(\1_n \otimes B) \left(- (\1_n\otimes U)\Theta_U(S(x)))\right) \vv^\flat(x)
\end{align}
and
\begin{align}\label{pair2}
\vv^\flat(x)&= \1_n\otimes C +(\1_n\otimes D)\frac{1}{\sqrt{2}}(\1+A(x))\vv(x) \notag\\
	&= \1_n\otimes C -(\1_n\otimes D) \left( (\1_n\otimes U)\Theta_U(S(x))\right)\vv^\flat(x).
\end{align}
Eliminate $\vv^\flat(x)$ from this pair of equations to obtain
\beq\label{11.2}
\ph(x)=a\1_n - (\1_n\otimes BU)S(x)(U)\left(\1 + (\1_n\otimes DU)S(x)(U)\right)\inv( \1_n\otimes C).
\eeq
Let
\beq\label{newp}
p=\bbm a& -BU\\ C& -DU\ebm = T \bbm 1&0\\0&-U \ebm \in\mathcal{L}(\C\oplus\ell^2),
\eeq
so that $p$ is a contraction on $\C\oplus\ell^2$. Equation \eqref{11.2} states that, for $x\in \bb$,
\begin{align}\label{11.5}
\ph(x)&= \mathcal{F}_{\1\otimes p}^u ( S(x)(U)) \notag \\
	&= \mathcal{F}_{\1\otimes p}^u \left( \Theta_U(S(x))\right).
\end{align}
 According to the definition  \eqref{defPhi}
\[
\Phi = \mathcal{F}_{\1\otimes p}^u \circ \Theta_U.
\]
Equation \eqref{11.5} states precisely that $\ph=\Phi\circ S$ on $\bb$.

It remains to check that $\Phi$ is an analytic nc-function on the  nc-domain $\Omega$ and is bounded by $1$.   It is clear that $\Phi$ is well defined and maps an element $g\in\Omega\cap \M_n^\infty$ into the closed unit ball of $\M_n$, so that $\Phi$ is graded. Clearly $\Phi$ is Fr\'echet differentiable on the open unit ball of $ A_n(\D)$ for each $n\geq 1$.  
By Proposition \ref{partnc}, for $g, h \in \Omega$,
\begin{align*}
\Phi(g\oplus h) &= \mathcal{F}_{\1\otimes p}^u((g\oplus h)(U)) = \mathcal{F}_{\1\otimes p}^u(g(U)\oplus h(U)) = \mathcal{F}_{\1\otimes p}^u(g(U))\oplus  \mathcal{F}_{\1\otimes p}^u(h(U)) \\
	&= \Phi(g) \oplus \Phi(h),
\end{align*}
and so $\Phi$ respects direct sums.  It also respects similarities.  Consider $g\in \Omega\cap A_n(\D)$ and an invertible matrix $s\in\M_n$ such that $s\inv gs\in\Omega$.  Note that, if $g=(g^0,g^1,g^2,\dots)\in \M_n$,
\[
\Theta_U (s\inv gs) = \sum_{j=0}^\infty (s\inv g^j s)\otimes U^j = (s\inv\otimes \1_{\ell^2}) \Theta_U(g) (s\otimes \1_{\ell^2}).
\]
Consequently
\begin{align*}
\Phi(s\inv gs) &= \mathcal{F}_{\1\otimes p}^u((s\inv g s)(U)) \\
	&= \mathcal{F}_{\1\otimes p}^u\left((s\inv\otimes \1_{\ell^2}) g(U) (s\otimes \1_{\ell^2})\right).
\end{align*}
Apply Proposition \ref{partnc} with $\mathcal{H}_1=\mathcal{K}_1=\C, \, \mathcal{H}_2=\mathcal{K}_2= \ell^2$ (recall equation \eqref{newp}) to obtain
\begin{align*}
\Phi(s\inv gs)	&= s\inv  \mathcal{F}_{\1\otimes p}^u(g(U)) s\\
	&= s\inv \Phi(g)s.
\end{align*}
Thus $\Phi$ is an nc-function on $\Omega$.
\end{proof}

In the course of the above proof the following realization formula was derived.
\begin{corollary}
For every symmetric function $\ph$ on $\bb$ bounded by $1$ in norm there exist a unitary operator $U$ on $\ell^2$ and a contraction $p$ on $\C\oplus \ell^2$ such that
\beq\label{realizph}
\ph =  \mathcal{F}_{\1\otimes p}^u \circ \Theta_U\circ S.
\eeq
\end{corollary}
This is just a restatement of equation \eqref{11.5}.  Diagrammatically, $U$ and $p$ satisfy

\begin{equation*}
\begin{array}{ccccccccc}
~ & ~ & \;\;\mathcal{M}^{\infty} ~ & ~ &~\\
~ & ~ &  \cup ~ & ~ &~\\
~ & ~ & \Omega  & ~ & ~ \\  
   ~ &  \vcenter{\llap{$\scriptstyle{S}$}}\nearrow & ~ & \searrow \vcenter{\rlap{$\scriptstyle{\Theta_U}$}}& ~ \\
  B^2 \;\; & ~ & ~ & 
~ & \;\; \bigcup_n \mathrm{ball} \;\mathcal{L}(\mathbb{C}^n \otimes \ell^2)\!\!\!  \\ 
~ &  \vcenter{\llap{$\scriptstyle{\varphi}$}}\searrow & ~ 
& \swarrow\vcenter{\rlap{$\scriptstyle{\mathcal{F}^u_{\1 \otimes p}}$}} & ~ \\ 
   ~ & ~ & \!\! \mathcal{M}^1\!\!  & ~ & ~ \\  
\end{array}
\end{equation*}
where $\mathrm{ball} \;\mathcal{L}(\mathbb{C}^n \otimes \ell^2)$ denotes the open unit ball of $\mathcal{L}(\C^n\otimes \ell^2)$.
\begin{remark} \rm  (1)  There is a trivial converse to Theorem \ref{realizn}.  If $\Phi:\Omega\to\M^1$ is a bounded analytic nc-function then $\Phi\circ S$ is a symmetric bounded analytic nc-function on $\bb$, with the same bound.\\

\noindent (2)  The realization formula \eqref{realizph} can be re-stated in terms of the {\em Redheffer product \cite{redheffer}}.  If $A,B$ are suitable $2\times 2$ operator matrices then  $B*A$ is the $2\times 2$ operator matrix with the property
\[
\mathcal{F}_{B*A}^u(X) = \mathcal{F}_{B}^u \circ \mathcal{F}_{A}^u(X)
\]
for every $X$ for which the expressions make sense.  In fact
\[
B*A = \bbm \mathcal{F}_B^u(A_{11}) & B_{12}(\1-A_{11}B_{22})\inv A_{12} \\
	A_{21}(1-B_{22}A_{11})\inv B_{21} & \mathcal{F}_A^\ell(B_{22}) \ebm.
\]
If we take
\[
A(x)=Q(x) \otimes \1_{\ell^2} = \bbm u&v\\v&u \ebm \otimes \1_{\ell^2}, \qquad
	B= \1_n\otimes p,
\]
then we find that, for $x\in\bb\cap\M_n^2$,
\begin{align*}
\ph(x) &= \mathcal{F}_B \circ \mathcal{F}_{A(x)} (\1_n\otimes U)\\
	&= \mathcal{F}_{B*A(x)}(\1_n \otimes U).
\end{align*}
Consequently
\beq\label{complicated}
\ph(x) = \mathcal{F}_{C(x)}(\1_n\otimes U)
\eeq
where
\begin{align*}\label{another}
C(x)&= B*A(x) \notag \\
	&=\bbm \mathcal{F}_{\1_n\otimes p}^u(u\otimes \1_{\ell^2}) & (\1_n\otimes p_{12})(\1-u\otimes p_{22})\inv (v\otimes \1_{\ell^2}) \\
(v\otimes \1_{\ell^2} )(1-u\otimes p_{22})\inv (\1_n\otimes p_{21}) & \mathcal{F}_{Q(x)\otimes \1_{\ell^2}}^\ell(\1_n\otimes p_{22})  \ebm
\end{align*}
and, as usual, $u=\half(x^1+x^2), \, v=\half(x^1-x^2)$.  The representation \eqref{complicated} differs from familiar realization formulae in that it is linear fractional not in $x$, but in $\1\otimes U$. \\ 

\noindent (3)  Since the operator $p$ in equation \eqref{defT} corresponds to  the Schur-class scalar function
\[
\psi(\la)=p_{11}+p_{12}\la(1-p_{22}\la)\inv p_{21}
\]
one might expect that $\Phi$ could be written in terms of $\psi$ and the functional calculus $\Theta_U$.  However, $\Phi$ depends on the particular realization of $\psi$; if
\[
q= (1\oplus s)\inv p (1\oplus s)
\]
for some invertible operator $s$ on $\ell^2$ then $\mathcal{F}_{\1\otimes q}^u \neq \mathcal{F}_{\1\otimes p}^u $ in general.
\end{remark}

As we observed in the introduction, $\Omega$ is not a true analogue of the symmetrised bidisc $\pi(\D^2)$ because the nc map $S:\bb\to\Omega$ is not surjective. To repair this failing we might replace $\Omega$ by its subset $S(\bb)$.  However, $S(\bb)$ is not an open subset of $\M^\infty$ in any natural topology.  

We ask: for a given bounded symmetric analytic nc-function $\ph$ on $\bb$,
is there a {\em unique}  analytic nc-function $\Phi:\Omega\to \M^1$ such that $\ph= \Phi\circ S$?

If one does not require $\Phi$ to be an nc-function then $\Phi$ is not unique.  Let $\ph$ be the zero function on $\bb$: then we may construct a non-zero analytic $\Phi$ on $\Omega$ such that $\Phi\circ S =\ph= 0$ as follows.
Fix $z_0\in\D, \, z_0 \neq 0$.
For $g\in \Omega\cap \M_n^\infty=\Omega\cap\M_n[[z]]$ let
\[
\Phi(g) = \left( \det(g(z_0)-g(0)) - z_0^n \frac{\det g'(0)}{\det(1-g(0)z_0)}\right) \1_n.
\]
$\Phi$ is well defined on $\Omega$ and is not identically zero at any level. It is easy to see that $\Phi(S(x))=0$ for all $x\in\bb$.  However $\Phi$ does not respect direct sums.

\begin{center} \bf Acknowledgements \end{center}
The first author was partially supported by National Science Foundation Grant on  Extending Hilbert Space Operators DMS 1068830.
The second author was partially supported by the UK Engineering and Physical Sciences Research Council grant  EP/K50340X/1 and London Mathematical Society grant 41219.


\begin{thebibliography}{1}

\bibitem{JJ} J. Agler and J. E. McCarthy, Global holomorphic functions in several non-commuting variables,  arXiv:1305.1636

\bibitem{AlpKal} D. Alpay and D. S. Kalyuzhnyi-Verbovetzkii, Matrix $J$-unitary non-commutative rational formal power series, in {\em The state space method: generalizations and applications}, {\bf OT 161} (2006) 49-113, Birkh\"auser, Basel.

\bibitem{BGM} J. A. Ball, G. Groenewald and T. Malakorn,  Conservative structured noncommutative multidimensional linear systems, in {\em The state space method: generalizations and applications}, {\bf OT 161} (2006) 179-223, Birkh\"auser, Basel.

\bibitem {gelfand} I. M. Gelfand, D. Krob, B. Leclerc, A. Lascoux, V. S. Retakh, J.-Y. Thibon,  Noncommutative
symmetric functions, {\em  Adv. Math.}, {\bf 112} (1995) 218--348. 

\bibitem{HelKlepMcC} J. W. Helton, I. Klep and S. McCullough, Proper analytic free maps, {\em J. Funct. Anal.} {\bf 260} (2011) 1476-1490.

\bibitem{KalVin} D. S. Kaliuzhnyi-Verbovetskyi and V. Vinnikov, {\em Foundations of Noncommutative Function Theory},   arXiv:1212.6345 .

\bibitem{popescu2006} G. Popescu, Free holomorphic functions on the unit ball of $B(\mathcal{H})^n$, {\em J. Funct. Anal.} {\bf 241} (2006) 268-333. 

\bibitem{redheffer} R. M. Redheffer, On a certain linear fractional transformation, {\em J.  Math. and Phys.}, {\bf 39} (1960) 269--286.

\bibitem{taylor} J. L. Taylor,  Functions of several non-commuting variables, {\em Bull. Amer. Math. Soc.} {\bf 79} (1973) 1--34.

\bibitem{voic} D. Voiculescu, Free analysis questions II: the Grassmannian completion and the series expansion at the origin, {\em J. Reine Angew. Math.} {\bf 645} (2010) 155-236.

\bibitem{wolf} M. Wolf,  Symmetric functions of noncommuting elements,  {\em Duke Math. J. } {\bf 2} (1936) 626--637.

\end{thebibliography}
\end{document}